\documentclass{amsart}
\usepackage{times,mathrsfs,array,amssymb,pb-diagram,epsfig}

\newcounter{lemma}
\newtheorem{Theorem}{Theorem}
\newtheorem{Lemma}[lemma]{Lemma}

\theoremstyle{definition}

\newtheorem{Remark}[lemma]{Remark}

\def\H{\mathbb H}
\def\O{\mathcal O}

\def\Q{\mathbb Q}
\def\R{\mathbb R}
\def\Z{\mathbb Z}
\def\Re{\mathrm{Re}\,}
\def\Im{\mathrm{Im}\,}
\def\FL#1{\left\lfloor #1\right\rfloor}

\def\mod{\  \mathrm{mod}\ }

\def\sgn{\mathrm{sgn\,}}
\def\Im{\mathrm{Im\,}}
\def\tr{\mathit{tr\,}}
\def\Nm{\text{N}}

\def\JS#1#2{\left(\frac{#1}{#2}\right)}
\def\SL{\mathrm{SL}}

\def\M#1#2#3#4{\begin{pmatrix}#1&#2\\#3&#4\end{pmatrix}}
\def\SM#1#2#3#4{\left(\begin{smallmatrix}#1&#2\\#3&#4\end{smallmatrix}
  \right)}
\def\wt{\widetilde}

\begin{document}

\title{Ramanujan-type identities for Shimura curves}

\author{Yifan Yang}
\address{Department of Applied Mathematics, National Chiao Tung
  University and National Center for Theoretical Sciences, Hsinchu,
  Taiwan 300}  
\email{yfyang@math.nctu.edu.tw}
\date{\today}
\thanks{The author was partially supported by Grant
  99-2115-M-009-011-MY3 of the National Science Council, Taiwan (R.O.C.).}

\begin{abstract} In 1914, Ramanujan gave a list of $17$ identities
  expressing $1/\pi$ as linear combinations of values of
  hypergeometric functions at certain rational numbers. Since then,
  identities of similar nature have been discovered by many authors.
  Nowadays, one of the standard approaches to this kind of identities
  uses the theory of modular curves. In this paper, we will consider
  the case of Shimura curves and obtain Ramanujan-type formulas
  involving special values of hypergeometric functions and products of
  Gamma values. These products of Gamma values are related to periods
  of elliptic curves with complex multiplication by $\Q(\sqrt{-3})$
  and $\Q(\sqrt{-4})$.
\end{abstract}
\maketitle

\begin{section}{Introduction}
Among all the amazing formulas discovered by Ramanujan, the identity
$$
  \sum_{n=0}^\infty\frac{(6n+1)(1/2)_n^3}{(n!)^3}\left(
  \frac14\right)^n=\frac4\pi
$$
involving hypergeometric functions and $\pi$ is perhaps one of the most
well-known, where $(a)_n$ is the Pochhammer symbol
$(a)_n=a(a+1)\ldots(a+n-1)$. In \cite{Ramanujan}, he gave $17$ such
identities. Since then, number theorists
\cite{Berndt,Borwein,Borwein2,Chan-Liaw,Chan-Liaw-Tan,Chan-Tanigawa,
Chan-Verrill,Chudnovsky,Zudilin} have produced many 
more identities of similar nature. (See \cite{survey} for a
comprehensive survey on the history of Ramanujan-type series for
$1/\pi$.)
%One of the most dramatic examples
%is
%$$
%  \sum_{n=0}^\infty\frac{(6n)!(13591409+545140134n)}{(3n)!(n!)^3(-640320)^{3n}}
% =\frac{426880\sqrt{10005}}\pi
%$$
%discovered by Chudnovsky and Chudnovsky \cite{}, each term of which
%yields $14$-digits of precision for $1/\pi$.
Nowadays, it is well-understood that Ramanujan's $1/\pi$-series are
related to the CM-theory and the Hecke theory of modular curves. Since
Shimura curves are generalizations of modular curves and there are
analogues of CM points and Hecke operators in the setting of
Shimura curves, one may wonder whether we have also Ramanujan-type
formulas for Shimura curves. In this paper,
we consider the Atkin-Lehner quotient $X_6^\ast$ of the Shimura curve
associated to a maximal order in the quaternion algebra of
discriminant $6$ over $\Q$ by all the Atkin-Lehner involutions. Using
the method developed in our earlier work \cite{Yang-Schwarzian} for
computing Hecke operators relative to the explicit bases of
automorphic forms in terms of solutions of Schwarzian differential
equations, we obtain Ramanujan-type identities for $X_6^\ast$.

To state the results, let us recall the Chowla-Selberg formula
\cite{Chowla-Selberg}, which states that if $E$ is an elliptic
curve over $\overline\Q$ with complex multiplication by an imaginary
quadratic number field $K$ of discriminant $d$, then, up to an algebraic
factor, the period of $E$ is
$$
  \Omega_d=\sqrt\pi\prod_{0<a<|d|}\Gamma\left(\frac a{|d|}\right)^{w_d\chi_d(a)/4h_d},
$$
where $\chi_d$ is the Kronecker character associated to $K$, $w_d$ is
the number of roots of unity in $K$, and $h_d$ is the class number of
$K$. Then we have the following Ramanujan-type identities for the
Shimura curve $X_6^\ast$.

\begin{Theorem} \label{theorem: identities}
Let
$$
  A_n=\frac{(1/12)_n(1/4)_n(5/12)_n}{(1/2)_n(3/4)_nn!}, \qquad
  A'_n=\frac{(7/12)_n(3/4)_n(11/12)_n}{(3/2)_n(5/4)_nn!},
$$
$$
  B_n=\frac{(1/12)_n(1/3)_n(7/12)_n}{(2/3)_n(5/6)_nn!}, \qquad
  B'_n=\frac{(5/12)_n(2/3)_n(11/12)_n}{(4/3)_n(7/6)_nn!},
$$
and
\begin{equation*}
\begin{split}
%  C_1&=(\sqrt3-\sqrt2)\frac{\Gamma(3/4)\Gamma(19/24)\Gamma(23/24)}
%  {\Gamma(5/4)\Gamma(13/24)\Gamma(17/24)}
%  =\frac4{\sqrt[4]{12}}\frac{\Gamma(3/4)^2}{\Gamma(1/4)^2}, \\
%  C_2&=(\sqrt2-1)\frac{\Gamma(5/6)\Gamma(17/24)\Gamma(23/24)}
%  {\Gamma(7/6)\Gamma(13/24)\Gamma(19/24)}
%  =\frac6{2^{7/6}}\frac{\Gamma(2/3)^3}{\Gamma(1/3)^3}.
  C_1=\frac{4}{\sqrt[4]{12}}\frac\pi{\Omega_{-4}^{2}}
     =\frac4{\sqrt[4]{12}}\frac{\Gamma(3/4)^2}{\Gamma(1/4)^2}, \qquad
  C_2=\frac{3}{\sqrt[6]2}\frac\pi{\Omega_{-3}^{2}}
     =\frac3{\sqrt[6]2}\frac{\Gamma(2/3)^3}{\Gamma(1/3)^3}.
\end{split}
\end{equation*}
Then
\begin{equation*}
\begin{split}
  \sum_{n=0}^\infty(R_1n+R_2)A_n\JS MN^n
&=R_3^{1/2}|M|^{3/4}N^{1/4}C_1 \\
  \sum_{n=0}^\infty(R_1n+R_1/2+R_2)A_n'\JS MN^n
&=R_3^{1/2}|M|^{1/4}N^{3/4}C_1^{-1}
\end{split}
\end{equation*}
hold for the following values of $M$, $N$, $R_1$, $R_2$ and $R_3$.
$$ \extrarowheight3pt
\begin{array}{r|rrrrr} \hline\hline
 D  & M & N & R_1 & R_2 & R_3 \\ \hline
-120 & -7^4 & 3^3\cdot5^3 & 74480 & 6860/3 & 5 \\
 -52 & 2^2\cdot3^7 & 5^6 & 64584 & 972 & 1\\
-132 & 2^4\cdot11^2 & 5^6 & 226512 & 1936 & 11 \\
 -75 & -11^4 & 2^{10}\cdot3^3\cdot5 & 946220 & 33275/3 & 2\\
 -43 & -3^7\cdot7^4 & 2^{10}\cdot5^6 & 159971868 & 3417309 & 2\\
 -88 & 3^7\cdot7^4 & 5^6\cdot11^3 & 324590112 & 3667356 & 11\\
-312 & 7^4\cdot23^4 & 5^6\cdot11^6 & 212826755232 & 901428696 & 6 \\
-148 & 2^2\cdot3^7\cdot7^4\cdot11^4 & 5^6\cdot17^6
     & 1027794972204 & 17528128002 & 1 \\
%-232 & -3^7\cdot7^4\cdot11^4\cdot19^4 & 5^6\cdot23^6\cdot29^3
%     & 2443514820645762864 & 37910578302227862 & 58\\
%-708 & 2^8\cdot7^4\cdot11^4\cdot47^4\cdot59^2 & 5^6\cdot17^6\cdot29^6
%     & 6167093737665689901072 & 99165724038088434496 & 59 \\
%-163 & -3^{11}\cdot7^4\cdot19^4\cdot23^4
%     &  2^{10}\cdot5^6\cdot11^6\cdot17^6
%     & 3071651418339975941652 & 15473113111724338791 & 2 \\
\hline\hline
\end{array}
$$
Also,
\begin{equation*}
\begin{split}
  \sum_{n=0}^\infty(R_1n+R_2)B_n\JS MN^n
&=R_3^{1/2}|M|^{2/3}N^{1/3}C_2 \\
  \sum_{n=0}^\infty(R_1n+R_1/3+R_2)B'_n\JS MN^n
&=R_3^{1/2}|M|^{1/3}N^{2/3}C_2^{-1}
\end{split}
\end{equation*}
hold for the following values of $M$, $N$, $R_1$, $R_2$, and
$R_3$.
$$ \extrarowheight3pt
\begin{array}{r|rrrrr} \hline\hline
 D & M & N & R_1 & R_2 & R_3 \\ \hline
-84 & 3^3 & 2^2\cdot7^2 & 4914 & 189/2 & 42 \\
-40 & -5^3 & 3^7 & 10200 & 175 & 2 \\
-51 & 2^{10} & 7^4 & 14688 & 384 & 2 \\
-19 & -2^{10} & 3^7 & 23712 & 896 & 2 \\
-168 & 5^6 & 7^2\cdot11^4 & 10773000 & 118125 & 42 \\
-228 & -3^6\cdot5^6 & 2^6\cdot7^4\cdot19^2
     & 3146100750 & 176023125/2 & 114 \\
-123 & 2^{10}\cdot5^6 & 7^4\cdot19^4 & 1323972000 & 18960000 & 2\\
-100 & -11^6 & 2^4\cdot3^7\cdot5\cdot7^4
     & 820122270 & 10907545/2 & 2\\
-147 & 2^{10}\cdot3^3\cdot5^6\cdot7 & 11^4\cdot23^4
     & 86680314000 & 2593080000 & 42 \\
-67 & -2^{16}\cdot5^6 & 3^7\cdot7^4\cdot11^4
    & 219324768000 & 2193920000 & 2 \\
-372 & 3^3\cdot5^6\cdot11^6 & 2^2\cdot7^4\cdot19^4\cdot31^2
     & 62650568301750 & 3597211344375/2 & 186 \\
-408 & 3^6\cdot5^6\cdot17^3 & 7^4\cdot11^4\cdot31^4
     & 33316016526000 & 159415290000 & 1\\
%-267 & 2^{16}\cdot5^6\cdot11^6 & 7^4\cdot31^4\cdot43^4
%     & 45769617921456000 & 1027334373120000 & 2\\
\hline\hline
\end{array}
$$
(Note that the entries in the first columns of the two tables indicate
the discriminants of the corresponding CM-points.)
\end{Theorem}

Note that the power series $\sum_n A_nx^n$ is in fact
$$
  \sum_{n=0}^\infty A_nx^n={}_3F_2\left(\frac1{12},\frac14,\frac5{12};
  \frac12,\frac34;x\right)={}_2F_1\left(\frac1{24},\frac5{24};
  \frac34;x\right)^2.
$$
The hypergeometric differential equation associated to
$_2F_1(1/24,5/24;3/4;x)$ is essentially the Schwarzian differential
equation associated to the Hauptmodul of $X_6^\ast$ that takes values
$0$, $1$, and $\infty$ at the CM-points of discriminants $-4$, $-24$,
and $-3$, respectively. Likewise, the power series $\sum_nB_nx^n$ is
the square of a $_2F_1$-hypergeometric function whose corresponding
hypergeometric differential equation is essentially the Schwarzian
differential equation associated to the Hauptmodul of $X_6^\ast$ that
has values $0$, $1$, and $\infty$ at the CM-points of discriminants
$-3$, $-24$, and $-4$, respectively. Thus, if we let
$\Omega_\infty=\infty$ be the period of a generalized elliptic curve
(which is natural in view of the well-known fact that the period of
the elliptic curve $y^2=x(x-1)(x-\lambda)$ is
$\pi{}_2F_1(1/2,1/2;1;\lambda)$), the general form of the Ramanujan's
identities and our identities is
$$
  \sum_{n=0}^\infty (R_1n+R_2) C_nx_0^n=\frac{R_3\pi}{\Omega_d^2},
$$
where $R_1,R_2,R_3\in\overline\Q$, $\sum_nC_nx^n$ is the power
series expansion of a meromorphic modular form of weight $2$ with
respect to a Hauptmodul $x$ of a modular curve or a Shimura curve such
that $x$ takes value $0$ at a CM-point of discriminant $d$ (possibly
$\infty$), and $x_0$ is the value of $x$ at some CM-point of
discriminant $d'\neq d$.

\begin{Remark} \label{remark: 1}
  There are three more pairs $(M,N)$ of integers such
  that the first set of equalities hold for some $R_1$, $R_2$, and
  $R_3$. However, a rigorous proof of these three cases will require
  the computation of Hecke operators $T_p$ on the spaces of
  automorphic forms of weight
  $8(p+1)$, for $p=29$, $31$, and $41$, respectively. Such a
  computation will take a considerable amount of time and computer
  resource. Numerically, we find that for
  $$ \extrarowheight3pt
  \begin{array}{r|rr} \hline\hline
  D & M & N \\ \hline
-232 & -3^7\cdot7^4\cdot11^4\cdot19^4 & 5^6\cdot23^6\cdot29^3 \\
-708 & 2^8\cdot7^4\cdot11^4\cdot47^4\cdot59^2 & 5^6\cdot17^6\cdot29^6 \\
-163 & -3^{11}\cdot7^4\cdot19^4\cdot23^4
     &  2^{10}\cdot5^6\cdot11^6\cdot17^6 \\
  \hline\hline
  \end{array}
  $$
  the first set of identities hold with
  $$ \extrarowheight3pt
  \begin{array}{r|rrr} \hline\hline
  D & R_1 & R_2 & R_3 \\ \hline
  -232 & 2443514820645762864 & 37910578302227862 & 58\\
  -708 & 6167093737665689901072 & 99165724038088434496 & 59 \\
  -163 & 3071651418339975941652 & 15473113111724338791 & 2 \\
  \hline\hline
  \end{array}
  $$
  Likewise, numerically, we find that the second set of identities
  holds with
  $M=2^{16}\cdot5^6\cdot11^6$, $N=7^4\cdot31^4\cdot43^4$,
  $R_1=45769617921456000$, $R_2=1027334373120000$, and $R_3=2$. This
  corresponds to the CM-point of discriminant $-267$. A rigorous proof
  of this case requires the computation of $T_{23}$ on the space of
  automorphic forms of weight $192$ on $X_6^\ast$.

  Note that if we allow $M/N$ to be irrational algebraic numbers, then
  there are infinitely many Ramanujan-type identities.
\end{Remark}

\begin{Remark} Observe that when a prime $p\neq 2,3$ divides $M$,
  the series in the theorem converges $p$-adically. One may
  wonder what kind of $p$-adic numbers they converge to. Quite
  interestingly and mysteriously, it appears that with $C_1$ and $C_2$
  replaced by certain products of $p$-adic Gamma values and possibly
  $R_3$ by different rational numbers, the identities also hold
  $p$-adically.

  Let $\Gamma_p(x)$ be the $p$-adic Gamma function defined by
  $$
    \Gamma_p(n)=(-1)^n\prod_{0<j<n,p\nmid j}j
  $$
  for positive integers $n$ and extended continuously to $\Z_p$. Let
  $$
    C_p=\frac{3^6}2\frac{\Gamma_p(2/3)^9}{\Gamma_p(1/3)^9}.
  $$
  Our numerical computation suggests that, for the second set of
  identities in Theorem \ref{theorem: identities}, when $p\neq 2,3$
  divides $M$,
  \begin{equation*}
  \begin{split}
    \sum_{n=0}^\infty(R_1n+R_2)B_n\left(\frac MN\right)^n
  &=\sqrt[6]{R_3^3M^4N^2C_p}, \\
    \sum_{n=0}^\infty(R_1n+R_1/3+R_2)B_n'\left(\frac MN\right)^n
  &=p\sqrt[6]{R_3^3M^2N^4/C_p}
  \end{split}
  \end{equation*}
  hold $p$-adically with the same $R_1$, $R_2$, but possibly different
  $R_3$. For example, for $p=5$, the equality holds with
  $$ \extrarowheight3pt
  \begin{array}{c|ccccccccc} \hline\hline
  D & -40 &-168&-228&-123&-147&-67&-372&-408&-267 \\ \hline
  R_3& 2 & 42 &-38&2 &21/2 &2 & -62 & -3 & 1/2 \\
  \hline\hline
  \end{array}
  $$
  (There is a unique sixth root such that the equality holds.) The
  equalities are verified up to $100$ $5$-adic digits. There also
  appears to be $p$-adic analogues for the first set of identities in
  Theorem \ref{theorem: identities} in which the limits are
  $\Gamma_p(3/4)/\Gamma_p(1/4)$ times some algebraic numbers. We
  believe that the identities should be related to $p$-adic periods
  of elliptic curves with CM, but we do not have a proof of these
  identities yet.

%  Note that there are no $p$-adic analogues in the cases of
%  classical modular curves because in those cases $M/N$ are always
%  reciprocals of algebraic integers and there are no primes $p$ that
%  make the series converge $p$-adically.
\end{Remark}

The plan for the rest of the paper is as follow. In Sections
\ref{subsection: X6}--\ref{subsection: Hecke}, we set up our
notations and review basic properties of the Shimura curve under
consideration. In particular, in Section \ref{subsection: space}, we
will give an explicit description of automorphic forms on the Shimura
curve in terms of hypergeometric functions. The computation of Hecke
operators relative to these automorphic forms is pivotal in the proof
of the Ramanujan-type identities. In Section \ref{subsection: main},
we state several general identities, the specialization of which gives
us the identities in Theorem \ref{theorem: identities}. The proof of
these general identities will be given in Section \ref{section:
  proof}. In Section \ref{section: proof 1}, we prove the identities
in Theorem \ref{theorem: identities}. The main task there is the
evaluation of certain automorphic functions on Shimura curves
associated to Eichler orders of the quaternion algebra of discriminant
$6$ over $\Q$ at CM-points. The evaluation relies on some auxiliary
polynomials, which will be listed in the appendix.
\end{section}

\begin{section}{Ramanujan-type formulas for Shimura curves}

\begin{subsection}{The Shimura curve $X_6^\ast$}
\label{subsection: X6}
Let the quaternion algebra of discriminant $6$ over $\Q$ be presented
by $B=\JS{-1,3}\Q$, i.e., the algebra generated by $I$ and $J$ over
$\Q$ with the relations
$$
  I^2=-1, \qquad J^2=3, \qquad IJ=-JI.
$$
Choose the embedding $\iota:B\to M(2,\R)$ to be
\begin{equation} \label{equation: iota}
  I\longmapsto\M0{-1}10, \qquad
  J\longmapsto\M{\sqrt3}00{-\sqrt3}
\end{equation}
as in \cite[Section 5.5.2]{Alsina-Bayer}. Fix the maximal order $\O$
to be $\Z+\Z I+\Z J+\Z(1+I+J+IJ)/2$. Then the image of the norm-one
group of $\O$ under the embedding $\iota$ is
$$
  \Gamma=\left\{\frac12\M\alpha\beta{-\beta'}{\alpha'}
  \in\SL(2,\R):~\alpha,\beta\in\Z[\sqrt3],~\alpha\equiv\beta\mod 2\right\},
$$
where $\alpha'$ and $\beta'$ denote the Galois conjugates of $\alpha$
and $\beta$, respectively. Let also
\begin{equation*}
\begin{split}
  \Gamma^\ast&=\left\{\frac1{\Nm(\gamma)^{1/2}}\iota(\gamma):
  \gamma\in N_B(\O),~\Nm(\gamma)>0\right\}
 =\gamma_1\Gamma\cup\gamma_2\Gamma\cup\gamma_3\Gamma\cup\gamma_6\Gamma,
\end{split}
\end{equation*}
where $N_B(\O)$ is the normalizer of $\O$ in $B$, and
\begin{equation} \label{equation: A-L matrices}
\begin{split}
  \gamma_1&=\iota(1)=\M1001, \qquad
  \gamma_2=\frac1{\sqrt2}\iota(1+I)=\frac1{\sqrt2}\M1{-1}11, \\
  \gamma_3&=\frac1{\sqrt3}\iota((3+3I+J+K)/2)
  =\frac1{2\sqrt3}\M{3+\sqrt3}{-3+\sqrt3}{3+\sqrt3}{3-\sqrt3}, \\
  \gamma_6&=\frac1{\sqrt6}\iota(3I+J)=\M0{-3+\sqrt3}{3+\sqrt3}0.
\end{split}
\end{equation}

Let $X_6=\Gamma\backslash\H$ and
$X_6^\ast=\Gamma^\ast\backslash\H$ be the Shimura curves associated
to $\O$ and $N_B(\O)$, respectively. In particular,
$X_6^\ast=X_6/W_6$ is the quotient of $X_6$ by all the
Atkin-Lehner involutions. A fundamental domain for $X^\ast_6$ is
given by

\centerline{\epsfig{file=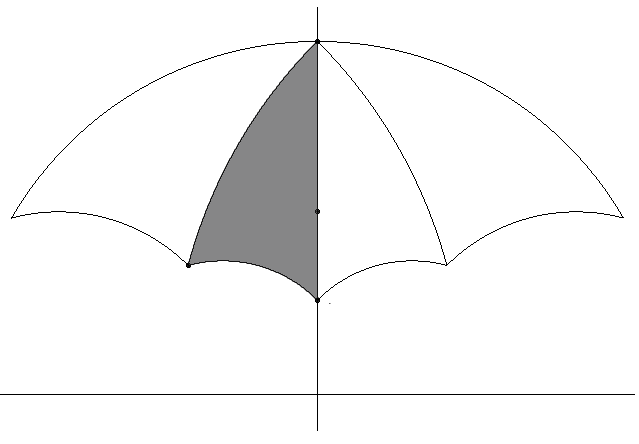,height=2.1in,width=3in}}

\noindent
(See \cite[Figure 5.1]{Alsina-Bayer} and \cite{Voight-computation}.)
Here the grey area represents a fundamental domain for $X_6^\ast$. The
grey area, together with the three white areas, forms a fundamental
domain for $X_6$. The four marked points on the boundary of the grey
area are
\begin{equation} \label{equation: elliptic points}
  P_4=i, \quad P_6=\frac{-1+i}{1+\sqrt3}, \quad (2-\sqrt3)i, \quad
  P_2=\frac{(\sqrt6-\sqrt2)i}2,
\end{equation}
respectively. The points $P_2$, $P_4$, and $P_6$ are representatives
of the elliptic points of orders $2$, $4$, and $6$, respectively.
They are CM-points of discriminants $-24$, $-4$, and $-3$,
respectively. Their isotropy subgroups are generated by
$$
  M_2=\frac1{\sqrt6}\M0{-3+\sqrt3}{3+\sqrt3}0, \qquad
  M_4=\frac1{\sqrt2}\M1{-1}11,
$$
and
\begin{equation*} % \label{equation: A6}
  M_6=\frac1{2\sqrt3}\M{3+\sqrt3}{3-\sqrt3}{-3-\sqrt3}{3-\sqrt3},
\end{equation*}
respectively. The point $(2-\sqrt 3)i$ is equivalent to $P_4$ through
$M_2$.
\end{subsection}

\begin{subsection}{Spaces of automorphic forms on $X_6^\ast$}
\label{subsection: space}

In \cite{Yang-Schwarzian}, we showed that when a Shimura curve has
genus $0$, all automorphic forms can be expressed in terms of
solutions of its Schwarzian differential equations. When a
Shimura curve corresponds a triangle group, this Schwarzian
differential equation is essentially a hypergeometric differential
equation. In this section, we review the properties relevant to our
consideration.

\begin{Lemma} \label{lemma: t'}
  Let $t(\tau)$ be the Hauptmodul of $X_6^\ast$ that
  takes values $0$, $1$, and $\infty$ at the elliptic points of order
  $4$, $2$, and $6$, respectively. Let
$$
  F_1(t)={}_2F_1\left(\frac1{24},\frac5{24};\frac34;t\right), \qquad
  F_2(t)=t^{1/4}{}_2F_1\left(\frac7{24},\frac{11}{24};\frac54;t\right),
$$
be two linearly independent solutions of
$$
  \theta\left(\theta-\frac14\right)-t\left(\theta+\frac1{24}\right)
  \left(\theta+\frac5{24}\right)=0, \qquad \theta=t\frac d{dt}.
$$
Then we have
\begin{equation} \label{equation: C}
  C\frac{F_2(t)}{F_1(t)}=\frac{\tau-P_4}{\tau-\overline P_4}
  =\frac{\tau-i}{\tau+i}
\end{equation}
where
\begin{equation*}
\begin{split}
 C&=\frac{P_2-P_4}{P_2-\overline P_4}
    \frac{\Gamma(3/4)\Gamma(19/24)\Gamma(23/24)}
    {\Gamma(5/4)\Gamma(13/24)\Gamma(17/24)} \\
  &=(\sqrt2-\sqrt3) \frac{\Gamma(3/4)\Gamma(19/24)\Gamma(23/24)}
    {\Gamma(5/4)\Gamma(13/24)\Gamma(17/24)}
   =-\frac4{\sqrt[4]{12}}\frac{\Gamma(3/4)^2}{\Gamma(1/4)^2}.
\end{split}
\end{equation*}
Moreover, we have
\begin{equation} \label{equation: t'}
  t'(\tau)=\frac{4t^{3/4}(1-t)^{1/2}}{C(P_4-\overline P_4)}
  (F_1(t)-CF_2(t))^2
 =\frac{2t^{3/4}(1-t)^{1/2}}{Ci}(F_1(t)-CF_2(t))^2,
\end{equation}
valid for $\tau$ in the fundamental domain such that $|t(\tau)|<1$.
Here $t(\tau)^{1/4}$ is defined in a way such that it becomes a
holomorphic function near $P_4$ and takes positive real values along
the boundary of the fundamental domain from $P_4$ to $P_2$. Likewise
$(1-t(\tau))^{1/2}$ is defined in a way such that it is a holomorphic
function near $P_2$ and takes positive values along the boundary from
$P_2$ and $P_6$.
\end{Lemma}

\begin{proof} The proof is similar to that of Lemma 14 in
  \cite{Yang-Schwarzian}, so we omit it. Here we just point
  out that the simplification of the constant $C$ uses the Gauss
  multiplication formula
  $$
    \Gamma(z)\Gamma(z+1/k)\ldots\Gamma(z+(k-1)/k)
   =(2\pi)^{(k-1)/2}k^{1/2-kz}\Gamma(kz)
  $$
  and the functional equation
  $$
    \Gamma(z)\Gamma(1-z)=\frac{\pi}{\sin(\pi z)},
  $$
  so that
  \begin{equation*}
  \begin{split}
    \left(\frac{c_{19}c_{23}}{c_{13}c_{17}}\right)^2
  &=\frac{c_3c_7c_{11}c_{15}c_{19}c_{23}}
    {c_1c_5c_9c_{13}c_{17}c_{21}}\cdot
    \frac{c_9c_{21}}{c_{3}c_{15}}\cdot
    \frac{c_1c_{23}}{c_7c_{17}}\cdot
    \frac{c_5c_{19}}{c_{11}c_{13}} \\
  &=6^{-1/2}\frac{\Gamma(3/4)}{\Gamma(1/4)}\cdot
    2^{-1/2}\frac{\Gamma(3/4)}{\Gamma(1/4)}\cdot
    \frac{\sin(\pi/24)}{\sin(7\pi/24)}\cdot
    \frac{\sin(5\pi/24)}{\sin(11\pi/24)} \\
  &=\frac{(\sqrt3-\sqrt2)^2}{\sqrt{12}}\frac{\Gamma(3/4)^2}{\Gamma(1/4)^2},
  \end{split}
  \end{equation*}
  where $c_j=\Gamma(j/24)$.
\end{proof}

\begin{Lemma}[{\cite[Theorem 9]{Yang-Schwarzian}}] \label{lemma: basis}
  Let $t(\tau)$ be the Hauptmodul of $X_6^\ast$ that
  takes values $0$, $1$, and $\infty$ at the elliptic points of order
  $4$, $2$, and $6$, respectively. For a positive even integer $k$,
  let
  $$
    d_k=1-k+\FL{\frac{k}4}+\FL{\frac{3k}8}+\FL{\frac{5k}{12}}
  $$
  be the dimension of the space $S_k(X_6^\ast)$ of automorphic
  forms of weight $k$ on $X_6^\ast$. Then a basis for
  $S_k(X_6^\ast)$ is
  $$
    t^jt^{\{3k/8\}}(1-t)^{\{k/4\}}(F_1(t)-CF_2(t))^k, \qquad
    j=0,\ldots,d_k-1.
  $$
\end{Lemma}
\end{subsection}

\begin{subsection}{CM-points on $X_6^\ast$}
For a negative integer $d$ with $d\equiv 0,1\mod 4$, we let $R_d$
denote the imaginary quadratic order of discriminant $d$ and $K_d$ be
the field of fractions of $R_d$. Recall that $K_d$ can be imbedded in
$B$ if and only if neither $2$ nor $3$ splits in $K_d$. An
\emph{optimal embedding} of $K_d$ into $B$ relative to $(\O,R_d)$ is
an embedding $\phi:K_d\hookrightarrow B$ such that
$\phi(K_d)\cap\O=\phi(R_d)$. We let $\Lambda(\O,d)$ denote the set of
all optimal embeddings relative to $(\O,R_d)$. Note that for each
$\phi\in\Lambda(\O,d)$, there is a unique point $\tau_\phi$ in the
upper half-plane such that
$$
  \iota(\phi(a))\tau_\phi=\tau_\phi
$$
for all $a\in K_d^\times$. We call $\tau_\phi$ the \emph{CM-point} of
discriminant $d$ associated to $\phi$. (In terms of moduli spaces, the
endomorphism ring of the corresponding abelian surface has center
$R_d$.)

Now it is clear that if $\phi\in\Lambda(\O,d)$, then so is
$\gamma^{-1}\phi\gamma$ for any $\gamma\in N_B(\O)$. In fact,
$\Lambda(\O,d)/N_B(\O)$ is a finite set. In the case when $d$ is a
fundamental discriminant, the cardinality is given by the formula
$$
  \left|\Lambda(\O,d)/N_B(\O)\right|=
  \begin{cases}
  1, &\text{if }d=-3,-4,-24, \\
  \displaystyle\frac{h(d)}4\left(1-\JS d2\right)\left(1-\JS d3\right),
  &\text{else}, \end{cases}
$$
where $h(d)$ denotes the class number of $R_d$. (See \cite{Vigneras}.)
From now on, we assume that the representatives of
$\Lambda(\O,d)/N_B(\O)$ are chosen such that their corresponding
CM-points are inside the fundamental domain drawn above. Moreover, if
$\tau$ is a CM-point of discriminant $d$ (fundamental or not) in the
fundamental domain with the corresponding optimal embedding
$\phi_\tau$, we let $\alpha_\tau$ denote
$$
  \alpha_\tau=\phi_\tau(\sqrt d).
$$

\begin{Lemma} \label{lemma: location}
  Let $t$ be the Hauptmodul of $X_6^\ast$ that takes
  values $0$, $1$, and $\infty$ at the elliptic points of orders $4$,
  $2$, and $6$, respectively. Suppose that $\tau$ is a CM-point of
  discriminant $d$ in the fundamental domain and write $\alpha_\tau$ as
  $a_0+a_1I+a_2J+a_3IJ$. If $t(\tau)$ takes a value in the line segment
  $[0,1]$, then $a_0=a_2=0$. If $t(\tau)$ takes a value in
  $[1,\infty)$, then $a_1=3a_3$. If $t(\tau)$ takes a negative value,
  then $a_0=0$ and $a_2=-a_3$.
\end{Lemma}

\begin{proof} It is clear that $\tr(\alpha_d)=0$. Thus, $a_0=0$.
  Also, the fixed point of $\iota(a_1I+a_2J+a_3IJ)$ in the upper
  half-plane is
  \begin{equation*} %\label{equation: tau0}
    \frac{a_2\sqrt3+\epsilon\sqrt d}{a_1+a_3\sqrt3},
  \end{equation*}
  where $\epsilon=\sgn(a_1+a_3\sqrt 3)$. Now if $t(\tau)$ takes a
  value in $[0,1]$, then $\tau$ lies on the vertical line segment from
  $P_4$ to $P_2$. Thus, $\Re\tau=0$ and we must have
  $a_2=0$.

  If $t(\tau)$ takes a value in $[1,\infty)$, then $\tau$ lies on the
  arc from $P_6$ to $P_2$, which is part of the circle of
  $|\tau|^2=2-\sqrt 3$. Thus,
  $$
    \frac{3a_2^2}{(a_1+a_3\sqrt3)^2}
   +\frac{a_1^2-3a_2^2-3a_3^2}{(a_1+a_3\sqrt3)^2}=2-\sqrt 3,
  $$
  which yields
  $$
    (a_1-3a_3)^2-(a_1-3a_3)(a_1-a_3)\sqrt 3=0.
  $$
  Therefore, $a_1=3a_3$.

  If $t(\tau)$ takes a negative value, then $\tau$ lies on the arc
  from $P_4$ to $P_6$, which is part of the circle
  $(\Re\tau-1)^2+(\Im\tau)^2=2$. Therefore,
  $$
    \left(\frac{a_2\sqrt3}{a_1+a_3\sqrt3}-1\right)^2
   +\frac{a_1^2-3a_2^2-3a_3^2}{(a_1+a_3\sqrt3)^2}=2.
  $$
  Simplifying, we find $a_2=-a_3$.
\end{proof}
\end{subsection}

\begin{subsection}{Hecke operators on the spaces of automorphic forms
    on $X_6^\ast$}
\label{subsection: Hecke}
To define the Hecke operator $T_p$ on $S_k(X_6^\ast)$ for a prime
$p$ relatively prime to $6$, we pick an element $\alpha$ in $\O$ of
norm $p$ and consider the double coset
\begin{equation} \label{equation: cosets}
  \Gamma_p^\ast=\Gamma^\ast\iota(\alpha)\Gamma^\ast
 =\gamma_1\Gamma_p\cup\gamma_2\Gamma_p\cup\gamma_3\Gamma_p\cup
  \gamma_6\Gamma_p,  
\end{equation}
where $\gamma_1,\gamma_2,\gamma_3,\gamma_6$ are defined by
\eqref{equation: A-L matrices}, and
$$
  \Gamma_p=\left\{\alpha\in\O: \Nm(\alpha)=p\right\}.
$$
Here $\Nm(\alpha)$ denotes the reduced norm of $\alpha$.
The group $\Gamma^\ast$ acts on the left of $\Gamma_p^\ast$ by
multiplication and the number of orbits is $p+1$.

For an automorphic form $G(\tau)$ of weight $k$ in $S_k(X_6^\ast)$
and an element $\gamma=\SM abcd\in\Gamma_p^\ast$, we define the slash
operator by
\begin{equation} \label{equation: slash}
  G(\tau)\Big|_k\gamma=\frac{(\det\gamma)^{k/2}}{(c\tau+d)^k}
  G\left(\frac{a\tau+b}{c\tau+d}\right).
\end{equation}
It is easy to show that if $\gamma_1$ and $\gamma_2$ is in the same
coset $\Gamma^\ast\backslash\Gamma_p^\ast$, then
$G|_k\gamma_1=G|_k\gamma_2$ and the function
$$
  \sum_{\gamma\in\Gamma^\ast\backslash\Gamma_p^\ast}G(\tau)\Big|_k\gamma
$$
is again an automorphic form of weight $k$ on $X_6^\ast$. In other
words, the mapping
$$
  T_p:G(\tau)\longmapsto p^{k/2-1}
  \sum_{\gamma\in\Gamma^\ast\backslash\Gamma_p^\ast}G(\tau)\Big|_k\gamma
$$
defines a linear transformation on $S_k(X_6^\ast)$. This linear
transformation is called the $p$th \emph{Hecke operator} of
$S_k(X_6^\ast)$. For general integers $n$ relatively prime to $6$, the
$n$th Hecke operator $T_n$ is defined similarly, with $\Gamma_p^\ast$
replaced by
$\Gamma_n^\ast=\Gamma_n\cup\gamma_2\Gamma_n\cup\gamma_3\Gamma_n
 \cup\gamma_6\Gamma_n$, where $\Gamma_n=\{\alpha\in\O: \Nm(\alpha)=n\}$.

Note that for an element $\gamma\in\gamma_e\Gamma_n$,
$e\in\{1,2,3,6\}$, the trace of $\sqrt e\gamma$ is an integer multiple
of $e$. Let $t=\tr(\sqrt e\gamma)/e$. If $t$ satisfies $(et)^2<4en$,
then the fixed point of $\gamma$ is a CM-point of discriminant
$(e^2t^2-4en)/f^2$ for some integer $f$.
\end{subsection}

\begin{subsection}{The main identities}
\label{subsection: main}
  Here we shall state the main identities, from which the identities
  in Theorem \ref{theorem: identities} are derived.

  Let the notations $B$, $\O$, and $X_6^\ast$ be defined as in Section
  \ref{subsection: X6} and $\iota$ be the embedding $B\hookrightarrow
  M(2,\R)$ given by \eqref{equation: iota}. Let the representatives
  $P_2$, $P_4$ and $P_6$ of elliptic points of orders $2$, $4$, $6$ be
  chosen as in \eqref{equation: elliptic points}. Let $t$ be the
  unique Hauptmodul of $X^\ast_6$ that takes values $0$, $1$, and
  $\infty$ at $P_4$, $P_2$, and $P_6$, respectively, and fix a nonzero
  automorphic form $F(\tau)$ of weight $8$ and a nonzero automorphic
  form $G(\tau)$ of weight $12$ on $X_6^\ast$.

  Let $D$ be the discriminant of an imaginary quadratic order $R_D$ with
  $K=\Q(\sqrt D)$ such that an optimal embedding $K\hookrightarrow B$
  relative to $(\O,R_D)$ exists. Let $\phi:K\hookrightarrow B$ be an
  optimal embedding relative to $(\O,R_D)$ such that the fixed point
  $\tau_0$ of $\iota(\phi(K^\ast))$ in the upper half-plane lies in
  the fundamental domain given in Section \ref{subsection: X6}. Write
  $\phi(\sqrt D)=a_1I+a_2J+a_3IJ$ and set
  $$
    \gamma=\iota(\phi(\sqrt D))
   =\M{a_2\sqrt 3}{-a_1+a_3\sqrt 3}{a_1+a_3\sqrt3}{-a_2\sqrt 3}.
  $$
  Changing $\phi$ to $-\phi$ if necessary, we assume that
  $a_1+a_3\sqrt 3>0$. Set
  $$
    \wt F(\tau)=F(\tau)\big|_8\gamma, \qquad
    \wt G(\tau)=G(\tau)\big|_{12}\gamma,
  $$
  where the slash operator is defined as in \eqref{equation: slash}.

  Now with the settings given as above, our main identities
  state as follows.

\begin{Theorem} \label{theorem: main identities}
  Set
  $$
%    C_1=(\sqrt3-\sqrt2)\frac{\Gamma(3/4)\Gamma(19/24)\Gamma(23/24)}
%   {\Gamma(5/4)\Gamma(13/24)\Gamma(17/24)}, \ \ 
%    C_2=(\sqrt2-1)\frac{\Gamma(5/6)\Gamma(17/24)\Gamma(23/24)}
%   {\Gamma(7/6)\Gamma(13/24)\Gamma(19/24)}.
    C_1=\frac4{\sqrt[4]{12}}\frac{\Gamma(3/4)^2}{\Gamma(1/4)^2},
    \qquad
    C_2=\frac3{\sqrt[6]2}\frac{\Gamma(2/3)^3}{\Gamma(1/3)^3}.
  $$
  Let $t_0=t(\tau_0)$. If $t_0$ is real and satisfies $0<t_0<1$, then
  \begin{equation} \label{equation: theorem 2 1}
    \sum_{n=0}^\infty\frac{(1/12)_n(1/4)_n(5/12)_n}{(1/2)_n(3/4)_nn!}
    \left(8Rt_0^{-1}n-RS\right)t_0^n=\frac{2a_3\sqrt3}{\sqrt{|D|}}C_1,
  \end{equation}
  and
  \begin{equation} \label{equation: theorem 2 2}
  \sum_{n=0}^\infty\frac{(7/12)_n(3/4)_n(11/12)_n}{(3/2)_n(5/4)_nn!}
  \left(8Rt_0^{-1}(n+1/2)-RS\right)t_0^n
 =\frac{2a_3\sqrt3}{\sqrt{|D|t_0}}C_1^{-1},
  \end{equation}
  where 
  \begin{equation} \label{equation: RS}
    R=|t_0^{3/4}(1-t_0)^{1/2}|, \qquad
    S=t'(\tau_0)^{-1}\frac{d}{d\tau}\log\frac{F(\tau)}{\wt F(\tau)}
    \Big|_{\tau=\tau_0}
     =\frac d{dt}\log\frac{F(\tau)}{\wt F(\tau)}\Big|_{t=t_0}.
  \end{equation}
  If $-1<t_0<0$, then
  \begin{equation} \label{equation: theorem 2 3}
    \sum_{n=0}^\infty\frac{(1/12)_n(1/4)_n(5/12)_n}{(1/2)_n(3/4)_nn!}
    \left(8R|t_0|^{-1}n+RS\right)t_0^n=\frac{2a_3\sqrt6}{\sqrt{|D|}}C_1.
  \end{equation}
  and
  \begin{equation} \label{equation: theorem 2 4}
    \sum_{n=0}^\infty\frac{(7/12)_n(3/4)_n(11/12)_n}{(3/2)_n(5/4)_nn!}
    \left(8R|t_0|^{-1}(n+1/2)+RS\right)t_0^n
   =\frac{2a_3\sqrt6}{\sqrt{|Dt_0|}}C_1^{-1},
  \end{equation}
  where $R$ and $S$ are defined in the same way as \eqref{equation:
    RS}.

  If $t_0$ is real and satisfies $|t_0|>1$, we set $s_0=1/t_0$. When
  $0<s_0<1$, we have
  \begin{equation} \label{equation: theorem 2 5}
    \sum_{n=0}^\infty\frac{(1/12)_n(1/3)_n(7/12)_n}{(2/3)_n(5/6)_nn!}
    \left(12R's_0^{-1}n-R'S'\right)s_0^n
   =\frac{2\sqrt6(a_2+a_3)}{\sqrt{|D|}}C_2,
  \end{equation}
  and
  \begin{equation} \label{equation: theorem 2 6}
    \sum_{n=0}^\infty\frac{(5/12)_n(2/3)_n(11/12)_n}{(4/3)_n(7/6)_nn!}
    \left(12R's_0^{-1}(n+1/2)-R'S'\right)s_0^n
   =\frac{2\sqrt 6(a_2+a_3)}{\sqrt{|D|s_0}}C_2^{-1},
  \end{equation}
  where
  \begin{equation} \label{equation: R'S'}
    R'=|s_0^{5/6}(1-s_0)^{1/2}|, \quad
    S'=s'(\tau_0)^{-1}\frac d{d\tau}\log\frac{G(\tau)}{\wt G(\tau)}
       \Big|_{\tau=\tau_0}=\frac d{ds}\log\frac{G(\tau)}{\wt G(\tau)}
    \Big|_{s=s_0}.
  \end{equation}
  When $-1<s_0<0$, we have
  \begin{equation} \label{equation: theorem 2 7}
    \sum_{n=0}^\infty\frac{(1/12)_n(1/3)_n(7/12)_n}{(2/3)_n(5/6)_nn!}
    \left(12R'|s_0|^{-1}n+R'S'\right)s_0^n
   =\frac{2\sqrt2(a_1-3a_3)}{\sqrt{|D|}}C_2,
  \end{equation}
  and
  \begin{equation} \label{equation: theorem 2 8}
    \sum_{n=0}^\infty\frac{(5/12)_n(2/3)_n(11/12)_n}{(4/3)_n(7/6)_nn!}
    \left(12R'|s_0|^{-1}(n+1/2)+R'S'\right)s_0^n
   =\frac{2\sqrt 2(a_1-3a_3)}{|D|^{1/2}|s_0|^{1/3}}C_2^{-1},
  \end{equation}
  where $R'$ and $S'$ are defined as \eqref{equation: R'S'}.
\end{Theorem}
We will prove the main identities in the next section.
\end{subsection}
\end{section}

\begin{section}{Proof of the main identities (Theorem \ref{theorem:
      main identities})}
\label{section: proof}
Let $D<0$ be the discriminant of an imaginary quadratic order $R_D$
such that an optimal embedding $\phi$ of relative to $(\O,R_D)$
exists. Let $\alpha=\phi(\sqrt{D})$ and $\gamma=\iota(\alpha)$ so that
the fixed point $\tau_0$ of $\gamma$ is a CM-point of discriminant $D$
on $X_6^\ast$. By conjugating by a suitable element of $\Gamma^\ast$,
we may assume that $\tau_0$ lies in the fundamental domain described
in Section \ref{subsection: X6}. Write $\gamma=\SM abcd$. Changing
$\phi$ to $-\phi$ if necessary, we assume that $c>0$.

Let $t$ be the Hauptmodul of $X_6^\ast$ that takes values $0$, $1$,
and $\infty$ at the elliptic points of orders $4$, $2$, and $6$,
respectively. Let
$$
  F_1(t)={}_2F_1\left(\frac1{24},\frac5{24};\frac34;t\right), \qquad
  F_2(t)=t^{1/4}{}_2F_1\left(\frac7{24},\frac{11}{24};\frac54;t\right),
$$
By Lemma \ref{lemma: basis}, the one-dimensional space $S_8(X_6^\ast)$
is spanned by
\begin{equation} \label{equation: main 0}
  F(\tau)=(F_1(t)-CF_2(t))^8, \qquad
%  C=(\sqrt2-\sqrt3) \frac{\Gamma(3/4)\Gamma(19/24)\Gamma(23/24)}
%    {\Gamma(5/4)\Gamma(13/24)\Gamma(17/24)}.
  C=-\frac4{\sqrt[4]{12}}\frac{\Gamma(3/4)^2}{\Gamma(1/4)^2}.
\end{equation}
(Note that $C=-C_1$.) Let
$$
  \wt F(\tau)=F(\tau)\big|_8\gamma=\frac{(\det\gamma)^4}{(c\tau+d)^8}
  F\left(\frac{a\tau+b}{c\tau+d}\right).
$$
Consider
$$
  \frac d{d\tau}\log\frac{F(\tau)}{\wt F(\tau)}
 =\frac{F'(\tau)}{F(\tau)}
 -\frac{\det\gamma}{(c\tau+d)^2}\frac{F'(\gamma\tau)}{F(\gamma\tau)}
 +\frac{8c}{c\tau+d}
$$
and evaluate the two sides at
$$
  \tau=\tau_0=\frac{a-d+\sqrt{(a-d)^2+4bc}}{2c}.
$$
We have
$$
  \det\gamma=|D|, \qquad
  c\tau_0+d=\frac{a+d+\sqrt{(a+d)^2-4\det\gamma}}2=\sqrt{D}.
$$
(Note that $a+d=\tr\gamma=0$.) Thus,
\begin{equation} \label{equation: main 1}
  \frac d{d\tau}\log\frac{F(\tau)}{\wt F(\tau)}\Bigg|_{\tau=\tau_0}
 =2\frac{F'(\tau_0)}{F(\tau_0)}+\frac{8c}{\sqrt{D}}.
\end{equation}
Set $t_0=t(\tau_0)$ and assume that $|t_0|<1$. For the right-hand
side of \eqref{equation: main 1}, we have, by \eqref{equation: main
  0},
\begin{equation} \label{equation: main RHS 1}
  \frac{F'(\tau_0)}{F(\tau_0)}=\frac8{F_1(t_0)-CF_2(t_0)}
  \frac d{d\tau}(F_1-CF_2)\Big|_{\tau=\tau_0}.
\end{equation}
Set
\begin{equation} \label{equation: u}
  u=\frac{\tau_0-P_4}{\tau_0-\overline P_4}=\frac{\tau_0-i}{\tau_0+i}.
\end{equation}
Then, by Lemma \ref{lemma: t'}, we have
\begin{equation} \label{equation: F1 F2}
  CF_2(t_0)=uF_1(t_0).
\end{equation}
Also, differentiating the two sides of \eqref{equation: C} in Lemma
\ref{lemma: t'} with respect to $\tau$ and then evaluating at
$\tau_0$, we obtain
\begin{equation*}
\begin{split}
  \frac{2i}{(\tau_0+i)^2}
&=\frac1{F_1(t_0)^2}\left(CF_1\frac d{d\tau}F_2-CF_2
  \frac d{d\tau}F_1\right)\Big|_{\tau=\tau_0} \\
&=\frac1{F_1(t_0)}\left(C\frac d{d\tau}F_2-u\frac d{d\tau}F_1\right)
  \Big|_{\tau=\tau_0}.
\end{split}
\end{equation*}
That is,
\begin{equation} \label{equation: dF2}
  C\frac d{d\tau}F_2\Big|_{\tau=\tau_0}
 =u\frac d{d\tau}F_1\Big|_{\tau=\tau_0}
  +\frac{2i}{(\tau_0+i)^2}F_1(t_0).
\end{equation}
Substituting this and \eqref{equation: F1 F2} into \eqref{equation:
  main RHS 1}, we get
\begin{equation} \label{equation: main RHS 2}
\begin{split}
  \frac{F'(\tau_0)}{F(\tau_0)}&=\frac8{(1-u)F_1(t_0)}\left(
  (1-u)\frac{dF_1}{d\tau}\Big|_{\tau=\tau_0}
 -\frac{2i}{(\tau_0+i)^2}F_1(t_0)\right) \\
&=\frac 8{F_1(t_0)}\frac{dF_1}{d\tau}\Big|_{\tau=\tau_0}
 -\frac8{\tau_0+i}.
\end{split}
\end{equation}
By \eqref{equation: t'} in Lemma \ref{lemma: t'}, we find
\begin{equation*}
\begin{split}
  \frac{dF_1}{d\tau}\Big|_{\tau=\tau_0}
&=\frac{dF_1(t)}{dt}\Big|_{t=t_0}t'(\tau_0)
 =\frac{2t_0^{3/4}(1-t_0)^{1/2}}{Ci}
  (F_1(t_0)-CF_2(t_0))^2\frac{dF_1(t)}{dt}\Big|_{t=t_0} \\
&=\frac{2R(1-u)^2}{Ci}F_1(t_0)^2
  \frac{dF_1(t)}{dt}\Big|_{t=t_0},
\end{split}
\end{equation*}
where
\begin{equation} \label{equation: A}
  R=t_0^{3/4}(1-t_0)^{1/2}.
\end{equation}
Plugging this into \eqref{equation: main RHS 2}, we arrive at
\begin{equation} \label{equation: main RHS 3}
  \frac{F'(\tau_0)}{F(\tau_0)}
 =\frac{16R(1-u)^2}{Ci}F_1(t_0)
  \frac{dF_1(t)}{dt}\Big|_{t=t_0}-\frac 8{\tau_0+i}.
\end{equation}

For the left-hand side of \eqref{equation: main 1},
we have, by \eqref{lemma: t'} in Lemma \ref{lemma: t'} again,
\begin{equation} \label{equation: main LHS 1}
\begin{split}
  \frac d{d\tau}\log\frac{F(\tau)}{\wt F(\tau)}\Big|_{\tau=\tau_0}
&=\frac d{dt}\log\frac{F(\tau)}{\wt F(\tau)}\Big|_{t=t_0}\cdot
  \frac {2R}{Ci}(F_1(t_0)-CF_2(t_0))^2
\end{split}
\end{equation}
Observe that $F(\tau)/\wt F(\tau)$ is an automorphic function on
$\Gamma^\ast\cap(\gamma^{-1}\Gamma^\ast\gamma)$, which is a subgroup
of finite index in $\Gamma^\ast$. Therefore, $F(\tau)/\wt F(\tau)$ is
an algebraic function of $t$, and so is
$d\log(F(\tau)/\wt F(\tau))/dt$. Let
\begin{equation} \label{equation: B}
  S=\frac d{dt}\log\frac{F(\tau)}{\wt F(\tau)}\Big|_{t=t_0}.
\end{equation}
Then \eqref{equation: main LHS 1} becomes
$$
  \frac d{d\tau}\log\frac{F(\tau)}{\wt F(\tau)}\Big|_{\tau=\tau_0}
 =\frac{2RS}{Ci}(1-u)^2F_1(t_0)^2.
$$
Combining this with \eqref{equation: main 1} and \eqref{equation: main
  RHS 3}, we obtain
\begin{equation} \label{equation: main 3}
  \frac{2RS}{Ci}(1-u)^2F_1(t_0)^2=
  \frac{32R}{Ci}(1-u)^2F_1(t_0)\frac{dF_1(t)}{dt}\Big|_{t=t_0}
 -\frac{16}{\tau_0+i}+\frac{8c}{\sqrt D}.
\end{equation}
Write $\phi(\sqrt D)=a_1I+a_2J+a_3IJ$ so that
$$
  \gamma=\M{a_2\sqrt3}{-a_1+a_3\sqrt3}{a_1+a_3\sqrt3}{-a_2\sqrt3}, \qquad
  \tau_0=\frac{a_2\sqrt3+\sqrt D}{a_1+a_3\sqrt 3}.
$$
Then
\begin{equation} \label{equation: (1-u)}
\begin{split}
&(1-u)^{-2}\left(\frac{16}{\tau_0+i}-\frac{8c}{\sqrt D}\right)
 =-\frac{(\tau_0+i)^2}4\left(\frac{16}{\tau_0+i}
  -\frac{8(a_1+a_3\sqrt3)}{\sqrt D}\right) \\
&\qquad\qquad=-2(\tau_0+i)\frac{a_1+a_3\sqrt3}{\sqrt D}
  \left(\frac{2\sqrt D}{a_1+a_3\sqrt3}
 -\frac{a_2\sqrt3+\sqrt D}{a_1+a_3\sqrt 3}-i\right) \\
&\qquad\qquad=\frac{2(a_1+a_3\sqrt3)}{\sqrt D}(\tau_0+i)
 (\overline\tau_0+i).
\end{split}
\end{equation}
Now it is necessary to consider the two cases $0<t_0<1$ and $-1<t_0<0$
separately.

Assume first $0<t_0<1$. By Lemma \ref{lemma: location}, we have
$a_2=0$ and
\begin{equation} \label{equation: tau0 1}
  \gamma=\M0{-a_1+a_3\sqrt3}{a_1+a_3\sqrt3}0, \qquad
  \tau_0=\frac{\sqrt D}{a_1+a_3\sqrt 3}.
\end{equation}
Thus, from \eqref{equation: (1-u)} we get
\begin{equation*}
\begin{split}
 (1-u)^{-2}\left(\frac{16}{\tau_0+i}-\frac{8c}{\sqrt D}\right)
 =\frac{2(a_1+a_3\sqrt3)}{\sqrt D}
  \left(\frac{|D|}{(a_1+a_3\sqrt3)^2}-1\right)
\end{split}
\end{equation*}
Now $\det\gamma=|D|$ implies that
$|D|=(a_1+a_3\sqrt3)(a_1-a_3\sqrt3)$. It follows that
$$
  (1-u)^{-2}\left(\frac{16}{\tau_0+i}-\frac{8c}{\sqrt D}\right)
 =\frac{2(a_1+a_3\sqrt3)}{\sqrt D}\left(
  \frac{a_1-a_3\sqrt3}{a_1+a_3\sqrt3}-1\right)
 =-\frac{4a_3\sqrt 3}{i\sqrt{|D|}}.
$$
Then we deduce from this and \eqref{equation: main 3} that
\begin{equation} \label{equation: main identity 1}
  16RF_1(t_0)\frac{dF_1(t)}{dt}\Big|_{t=t_0}
  -RSF_1(t_0)^2=-\frac{2a_3\sqrt3}{\sqrt{|D|}}C,
\end{equation}
where $R$ and $S$ are the numbers defined in \eqref{equation: A} and
\eqref{equation: B}, respectively.

Furthermore, from \eqref{equation:
  t'} and \eqref{equation: dF2}, we have
\begin{equation*}
\begin{split}
  C\frac{dF_2(t)}{dt}\Big|_{t=t_0}&=u\frac{dF_1(t)}{dt}\Big|_{t=t_0}
  +\frac{2iF_1(t_0)}{(\tau_0+i)^2}\frac{Ci}{2R(F_1(t_0)-CF_2(t_0))^2}
  \\
&=u\frac{dF_1(t)}{dt}\Big|_{t=t_0}+\frac{2i}{(\tau_0+i)^2}
  \frac{Ci}{2R(1-u)^2F_1(t_0)} \\
&=u\frac{dF_1(t)}{dt}\Big|_{t=t_0}+\frac C{4RF_1(t_0)}.
\end{split}
\end{equation*}
From this, \eqref{equation: F1 F2}, and
\eqref{equation: main identity 1}, it follows that 
\begin{equation} \label{equation: main tmp}
  16R\frac CuF_2(t_0)\left(\frac Cu\frac{dF_2(t)}{dt}\Big|_{t=t_0}
  -\frac C{4RuF_1(t_0)}\right)-RS\frac{C^2}{u^2}F_2(t_0)^2
 =-\frac{2a_3\sqrt3}{\sqrt{|D|}}C.
\end{equation}
Equivalently, we have
\begin{equation} \label{equation: F2 identity 1}
  16RF_2(t_0)\frac{dF_2(t)}{dt}\Big|_{t=t_0}-RSF_2(t_0)^2
 =\left(4u-\frac{2u^2a_3\sqrt3}{\sqrt{|D|}}\right)C^{-1}.
\end{equation}
Now from \eqref{equation: u}, \eqref{equation: tau0 1}, and the fact
$|D|=a_1^2-3a_3^2$, we have
\begin{equation*}
\begin{split}
  u&=\frac{\sqrt{|D|}-a_1-a_3\sqrt3}{\sqrt{|D|}+a_1+a_3\sqrt3}
    =\frac{(\sqrt{|D|}-a_1-a_3\sqrt3)^2}{|D|-(a_1+a_3\sqrt3)^2} \\
   &=\frac{a_1^2-3a_3^2-2(a_1+a_3\sqrt3)\sqrt{|D|}+(a_1+a_3\sqrt3)^2}
     {-2a_1a_3\sqrt3-6a_3^2}=-\frac{a_1-\sqrt{|D|}}{a_3\sqrt3}
\end{split}
\end{equation*}
so that
\begin{equation*}
\begin{split}
  4u-\frac{2u^2a_3\sqrt3}{\sqrt{|D|}}
&=\frac{a_1-\sqrt{|D|}}{a_3\sqrt3}
  \left(-4-2\frac{a_1-\sqrt{|D|}}{\sqrt{|D|}}\right) \\
&=-2\frac{(a_1-\sqrt{|D|})(a_1+\sqrt{|D|})}{a_3\sqrt{3|D|}}
 =-\frac{2a_3\sqrt3}{\sqrt{|D|}}.
\end{split}
\end{equation*}
Substituting this into \eqref{equation: F2 identity 1}, we get
\begin{equation} \label{equation: main identity 2}
  16RF_2(t_0)\frac{dF_2(t)}{dt}\Big|_{t=t_0}-RSF_2(t_0)^2
 =-\frac{2a_3\sqrt3}{C\sqrt{|D|}}.
\end{equation}

Now we consider the cases $-1<t_0<0$. By Lemma \ref{lemma: location},
we have $a_2=-a_3$ so that
$$
  \gamma=\M{-a_3\sqrt3}{-a_1+a_3\sqrt3}{a_1+a_3\sqrt3}{a_3\sqrt3},
  \qquad \tau_0=\frac{-a_3\sqrt3+\sqrt{D}}{a_1+a_3\sqrt3}.
$$
In this case, we have $|D|=a_1^2-6a_3^2$ and \eqref{equation: (1-u)}
becomes
\begin{equation*}
\begin{split}
 (1-u)^{-2}\left(\frac{16}{\tau_0+i}-\frac{8c}{\sqrt D}\right)
&=\frac{2(a_1+a_3\sqrt3)}{\sqrt D}\left(
  \frac{a_1^2-3a_3^2}{(a_1+a_3\sqrt3)^2}
 -\frac{2a_3\sqrt3 i}{a_1+a_3\sqrt3}-1\right) \\
&=-\frac{4a_3\sqrt3}{\sqrt D}(1+i)
 =-\frac{4a_3\sqrt 6}{i\sqrt{|D|}}e^{2\pi i/8}.
\end{split}
\end{equation*}
Instead of \eqref{equation: main identity 1}, we now have
\begin{equation} \label{equation: main identity 1'}
  16RF_1(t_0)\frac{dF_1(t)}{dt}\Big|_{t=t_0}
  -RSF_1(t_0)^2=-\frac{2a_3\sqrt6}{\sqrt{|D|}}e^{2\pi i/8}C.
\end{equation}
Likewise, instead of \eqref{equation: main tmp}, we have
$$
  16R\frac CuF_2(t_0)\left(\frac Cu\frac{dF_2(t)}{dt}\Big|_{t=t_0}
  -\frac C{4RuF_1(t_0)}\right)-RS\frac{C^2}{u^2}F_2(t_0)^2
 =-\frac{2a_3\sqrt6}{\sqrt{|D|}}e^{2\pi i/8}C,
$$
or equivalently,
\begin{equation} \label{equation: tmp 2}
  16RF_2(t_0)\frac{dF_2(t)}{dt}\Big|_{t=t_0}-RSF_2(t_0)^2
 =\left(4u-\frac{2u^2a_3\sqrt6}{\sqrt{|D|}}e^{2\pi ii/8}\right)C^{-1}.
\end{equation}
Now
\begin{equation*}
\begin{split}
  u&=\frac{-a_3\sqrt3+i(\sqrt{|D|}-a_1-a_3\sqrt3)}
     {-a_3\sqrt3+i(\sqrt{|D|}+a_1+a_3\sqrt3)}
  =\frac{(a_3\sqrt3+i(a_1+a_3\sqrt3))^2+|D|}
  {3a_3^2+(\sqrt{|D|}+a_1+a_3\sqrt3)^2} \\
&=\frac{-2a_1a_3\sqrt3-6a_3^2+2a_3\sqrt3(a_1+a_3\sqrt3)i}
  {2a_1^2+2a_1a_3\sqrt3+2(a_1+a_3\sqrt3)\sqrt{|D|}} \\
&=\frac{a_3\sqrt3(-1+i)}{a_1+\sqrt{|D|}}
 =-\frac{a_1-\sqrt{|D|}}{a_3\sqrt6}e^{-2\pi i/8}
\end{split}
\end{equation*}
and hence
\begin{equation*}
\begin{split}
  4u-\frac{2u^2a_3\sqrt6}{\sqrt{|D|}}e^{2\pi i/8}
&=2u\left(2+\frac{a_1-\sqrt{|D|}}{\sqrt{|D|}}\right) \\
&=-2e^{-2\pi i/8}\frac{a_1-\sqrt{|D|}}{a_3\sqrt6}\frac{a_1+\sqrt{|D|}}
  {\sqrt{|D|}}
 =-2e^{-2\pi i/8}\frac{a_3\sqrt6}{\sqrt{|D|}}
\end{split}
\end{equation*}
Then \eqref{equation: tmp 2} becomes
\begin{equation} \label{equation: main identity 2'}
  16RF_2(t_0)\frac{dF_2(t)}{dt}\Big|_{t=t_0}-RSF_2(t_0)^2
 =-\frac{2a_3\sqrt6}{\sqrt{|D|}}e^{-2\pi i/8}C^{-1}.
\end{equation}

Now by Clausen's identity \cite{SpecialFunctions}
\begin{equation} \label{equation: Clausen}
  _2F_1(\alpha,\beta;\alpha+\beta+1/2;z)^2
 ={}_3F_2(2\alpha,\alpha+\beta,2\beta;2\alpha+2\beta,
  \alpha+\beta+1/2;z),
\end{equation}
we may write
\begin{equation*}
\begin{split}
  F_1(t_0)^2&={}_2F_1\left(\frac1{24},\frac5{24};\frac34;t_0\right)^2
 ={}_3F_2\left(\frac1{12},\frac14,\frac5{12};
  \frac12,\frac34;t_0\right) \\
&=\sum_{n=0}^\infty\frac{(1/12)_n(1/4)_n(5/12)_n}{(1/2)_n(3/4)_nn!}
  t_0^n,
\end{split}
\end{equation*}
and
\begin{equation*}
\begin{split}
  2F_1(t_0)\frac{dF_1(t)}{dt}\Big|_{t=t_0}
 =\sum_{n=0}^\infty\frac{(1/12)_n(1/4)_n(5/12)_n}{(1/2)_n(3/4)_nn!}
  nt_0^{n-1}.
\end{split}
\end{equation*}
Substituting these two expressions into \eqref{equation: main identity
  1}, we find that for the cases $0<t_0<1$,
\begin{equation*}
\begin{split}
  \sum_{n=0}^\infty\frac{(1/12)_n(1/4)_n(5/12)_n}{(1/2)_n(3/4)_nn!}
  \left(8Rt_0^{-1}n-RS\right)t_0^n=-\frac{2a_3\sqrt3}{\sqrt{|D|}}C.
\end{split}
\end{equation*}
This proves \eqref{equation: theorem 2 1}. For the cases $-1<t_0<0$,
we note that, according to the description of $t^{1/4}$ given in Lemma
\ref{lemma: t'}, we have $R=-e^{2\pi i/8}|R|$. In these cases,
\eqref{equation: main identity 1'} becomes
$$
  \sum_{n=0}^\infty\frac{(1/12)_n(1/4)_n(5/12)_n}{(1/2)_n(3/4)_nn!}
  \left(8|R||t_0|^{-1}n+|R|S\right)t_0^n=-\frac{2a_3\sqrt6}{\sqrt{|D|}}C.
$$
This is \eqref{equation: theorem 2 3} in the theorem.

Likewise, by \eqref{equation: Clausen} again, we have
\begin{equation*}
\begin{split}
  F_2(t_0)^2&=t_0^{1/2}
  {}_2F_1\left(\frac7{24},\frac{11}{24};\frac54;t_0\right)^2
 =t_0^{1/2}{}_3F_2\left(\frac7{12},\frac34,\frac{11}{12};
  \frac32;\frac54;t_0\right)\\
&=t_0^{1/2}\sum_{n=0}^\infty
  \frac{(7/12)_n(3/4)_n(11/12)_n}{(3/2)_n(5/4)_nn!}t_0^n
\end{split}
\end{equation*}
and
$$
  2F_2(t_0)\frac{dF_2(t)}{dt}\Big|_{t=t_0}
 =t_0^{1/2}\sum_{n=0}^\infty
  \frac{(7/12)_n(3/4)_n(11/12)_n}{(3/2)_n(5/4)_nn!}(n+1/2)t_0^{n-1}
$$
Substituting these two into \eqref{equation: main identity 2}, we
obtain, for the cases $0<t_0<1$,
$$
  \sum_{n=0}^\infty\frac{(7/12)_n(3/4)_n(11/12)_n}{(3/2)_n(5/4)_nn!}
  \left(8Rt_0^{-1}(n+1/2)-RS\right)t_0^n
 =-\frac{2a_3\sqrt3}{\sqrt{|D|t_0}}C^{-1}.
$$
This proves \eqref{equation: theorem 2 2}. For the cases $-1<t_0<0$,
we have $t_0^{1/2}=e^{-2\pi i/4}|t_0|$ and $R=-e^{2\pi i/8}|R|$, and
\eqref{equation: main identity 2'} yields
$$
  \sum_{n=0}^\infty\frac{(7/12)_n(3/4)_n(11/12)_n}{(3/2)_n(5/4)_nn!}
  \left(8|R||t_0|^{-1}(n+1/2)+|R|S\right)t_0^n
 =-\frac{2a_3\sqrt6}{\sqrt{|Dt_0|}}C^{-1}.
$$
This is \eqref{equation: theorem 2 4} in the theorem.

The proof of the identities
\eqref{equation: theorem 2 5}--\eqref{equation: theorem 2 8} is very
similar to the proof of
\eqref{equation: theorem 2 1}--\eqref{equation: theorem 2 4} and is
skipped.
%Here we
%only mention the main difference and omit the rest of details. Instead
%of Lemma \ref{lemma: t'}, we now use the following properties.

%Let $s(\tau)$ be the Hauptmodul of $X_6^\ast$ that takes values $0$,
%$1$, and $\infty$ at the elliptic points of order $6$, $2$, and $4$,
%respectively. Let
%$$
%  G_1(s)={}_2F_1\left(\frac1{24},\frac7{24};\frac56;s\right), \qquad
%  G_2(s)={}_2F_1\left(\frac5{24},\frac{11}{24};\frac76;s\right).
%$$
%Then we have
%$$
%  C\frac{G_2(s)}{G_1(s)}=\frac{\tau-P_6}{\tau-\overline P_6},
%$$
%where
%\begin{equation*}
%\begin{split}
%  C'&=\frac{P_2-P_6}{P_2-\overline P_6}
%     \frac{\Gamma(5/6)\Gamma(17/24)\Gamma(23/24)}
%          {\Gamma(7/6)\Gamma(13/24)\Gamma(19/24)} \\
%   &=e^{-2\pi i/8}(\sqrt 2-1)\frac{\Gamma(5/6)\Gamma(17/24)\Gamma(23/24)}
%          {\Gamma(7/6)\Gamma(13/24)\Gamma(19/24)}.
%\end{split}
%\end{equation*}
%Then the one-dimensional space $S_{12}(X_6^\ast)$ is spanned by
%$$
%  G(\tau)=(G_1(s)-C'G_2(s))^{12}.
%$$
%Moreover, we have
%$$
%  s'(\tau)=\frac{6s^{5/6}(1-s)^{1/2}}{C'(P_6-\overline P_6)}
%  (G_1(s)-C'G_2(s))^2,
%$$
%valid for $\tau$ in the fundamental domain such that $|s(\tau)|<1$.
\end{section}

\begin{section}{Proof of Theorem \ref{theorem: identities}}
\label{section: proof 1}
Let $t(\tau)$ be the Hauptmodul of $X_6^\ast$ that takes values $0$,
$1$, and $\infty$ at the elliptic points of order $4$, $2$, and $6$,
respectively. There are finitely many discriminants $D$ such that the
number of CM-point of discriminant $D$ on $X_6^\ast$ is one. These
discriminants are given in the first columns of the tables in Theorem
\ref{theorem: identities} and Remark \ref{remark: 1}. The values of
$t$ at these CM-points were determined numerically by Elkies
\cite{Elkies} and later proved rigorously by Errthum \cite{Errthum}
using Borcherds forms. These values are given by $M/N$, where $M$ and
$N$ are the integers from the tables. Now according to Theorem
\ref{theorem: main identities}, to prove the identity associated to
the discriminant $D$, the main task is the evaluation of the
constant $S$ in \eqref{equation: RS} or $S'$ in \eqref{equation:
  R'S'}. Here we will work out two cases and omit the others.

Let $D=-120$ and choose the optimal embedding $\phi:\Q(\sqrt{-30})\to
B$ relative to $(\O,\Z[\sqrt{-30}])$ to be the one determined by
$\phi(\sqrt{-30})=6I-J+IJ$. Let
$$
  \gamma=\iota(\phi(\sqrt{-120}))
 =2\M{-\sqrt3}{-6+\sqrt3}{6+\sqrt3}{\sqrt3},
$$
and $\tau_0=(-\sqrt3+\sqrt{-30})/(6+\sqrt3)$ be the fixed point of
$\gamma$ in the upper half-plane. Then $\tau_0$ is the representative
of the CM-point of discriminant $-120$ in the fundamental domain
described in Section \ref{subsection: X6}. According to Theorem
\ref{theorem: main identities}, we need to evaluate the function
$$
  f(\tau)=t'(\tau)^{-1}\frac d{d\tau}\log\frac{F(\tau)}{\wt F(\tau)}
  =\frac d{dt}\log\frac{F(\tau)}{\wt F(\tau)}
$$
at $\tau_0$, where $F(\tau)$ is any nonzero automorphic form of weight
$8$ on $X_6^\ast$ and $\wt F(\tau)=F(\tau)\big|_8\gamma$. In order to
do so, we note that $g(\tau)$ is an automorphic function on the
subgroup $\Gamma^\ast\cap(\gamma^{-1}\Gamma^\ast\gamma)$ and hence
$t(\tau)$ and $g(\tau)$ satisfy a polynomial relation $P(t,f)=0$ for
some polynomial $P(x,z)$. Then $f(\tau_0)$ will be a root of the
polynomial $P(t(\tau_0),z)$ in $z$.

Here, to determine the polynomial $P(x,y)$, we observe that the matrix
$\gamma_0=\gamma/\sqrt6$ lies in $\Gamma_5^\ast$, where
$\Gamma_p^\ast$ is defined by \eqref{equation: cosets}, and can be
taken to be one of the coset representatives in
$\Gamma^\ast\backslash\Gamma_5^\ast$ defining the $5$th Hecke operator
$T_5$ on $X_6^\ast$. Let $\gamma_1,\ldots,\gamma_5$ be the other coset
representative defining $T_5$. Then any symmetric sum of
$F(\tau)\big|_8\gamma_j$ will be an automorphic function on $X_6^\ast$
and hence equal to a rational function in $t$. In particular, there is
a polynomial $Q(x,y)$ of degree $6$ in $y$ such that $Q(t, F/\wt F)=0$.
Then we have
\begin{equation*}
\begin{split}
  0&=\frac d{d\tau}Q(t(\tau),F(\tau)/\wt F(\tau)) \\
   &=t'(\tau)\frac\partial{\partial x}Q(x,y)
     \Big|_{x=t(\tau),y=F(\tau)/\wt F(\tau)}
    +\left(\frac{F(\tau)}{\wt F(\tau)}\right)'
     \frac\partial{\partial y}Q(x,y)
     \Big|_{x=t(\tau),y=F(\tau)/\wt F(\tau)}.
\end{split}
\end{equation*}
Eliminating the variable $y$ from the two polynomials
\begin{equation} \label{equation: Q}
  Q(x,y), \quad \frac{\partial}{\partial x}Q(x,y)
  +yz\frac\partial{\partial y}Q(x,y),
\end{equation}
we get the polynomial $P(x,z)$ satisfying $P(t,f)=0$.

From the above discussion, we see that the problem of proving the
identity for the discriminant $-120$ boils down to the determination
of the polynomial $Q(x,y)$ satisfying $Q(t,F/\wt F)=0$. This is where
the method of computing Hecke operators developed in
\cite{Yang-Schwarzian} comes in.

Let $u=-540/t$ be the Hauptmodul of $X_6^\ast$ that takes values $0$,
$-540$, and $\infty$ at $P_6$, $P_2$, and $P_4$, respectively. By
Equation (16) of \cite{Yang-Schwarzian}, for a positive even integer
$k$, a basis for the space of automorphic forms of weight $k$ on
$X_6^\ast$ is given by
\begin{equation*}
\begin{split}
  g_{k,\ell}&=u^{\{5k/12\}}(1+u/540)^{\{k/4\}}u^\ell \\
  &\qquad\times\left({}_2F_1\left(\frac1{24},\frac7{24};\frac56;
   -\frac u{540}\right)-C'u^{1/6}
   {}_2F_1\left(\frac5{24},\frac{11}{24};\frac76;-\frac u{540}\right)
   \right)^k,
\end{split}
\end{equation*}
$\ell=0,\ldots,d_k-1=\FL{5k/12}+\FL{l/4}+\FL{3k/8}-k$. In Section 4 of
\cite{Yang-Schwarzian}, we discussed how to compute Hecke 
operators with respect to these bases. The results relevant to our
problem here are given in the following table.
$$ \extrarowheight3pt
\begin{array}{c|l} \hline\hline
k & M \\ \hline
8 & -114 \\
16& 77646 \\
24& \displaystyle\begin{pmatrix}
       10980750 &  3111696/5 \\
    55987200000 &   14267406 \end{pmatrix} \\
32& \displaystyle\begin{pmatrix}
    105068988750 &  376515216/5 \\
  12317184000000 & -39127734834 \end{pmatrix} \\
40& \displaystyle\begin{pmatrix}
    -70619784011250  &  45558341136/5 \\
  30930128640000000  &  36422537206926 \end{pmatrix}\\
48& \displaystyle\begin{pmatrix}
       100480725468750 &       225950546273760  &      2420662999104/5 \\
    512317872000000000 &     22159766272716750  &      5512559277456/5 \\
2612138803200000000000 &  -7950573190656000000  &   -23013714467131314
  \end{pmatrix} \\
\hline\hline
\end{array}
$$
That is, we have
$$
  T_5\begin{pmatrix}g_{k,0}\\ \vdots\\ g_{k,d_k-1}\end{pmatrix}
=M\begin{pmatrix}g_{k,0}\\ \vdots\\ g_{k,d_k-1}\end{pmatrix}.
$$
Now we choose the nonzero automorphic form $F$ of weight $8$ to be
$g_{8,0}$. From the description of $g_{k,\ell}$, we find
$$
  F^{k/8}=g_{k,d_k-1}
$$
for any weight $k$ that is a multiple of $8$. Then the table above
gives us
$$
  T_5F=-114F, \qquad T_5F^2=77646F^2
$$
\begin{equation*}
\begin{split}
  T_5F^3&=(55987200000/u+14267406)F^3=(-103680000t+14267406)F^3, \\
  T_5F^4&=(-22809600000t-39127734834)F^4, \\
  T_5F^5&=(-57278016000000t+36422537206926)F^5, \\
  T_5F^6&=(8957952000000000t^2+14723283686400000t-23013714467131314)F^6.
\end{split}
\end{equation*}
By the definition of Hecke operators, we have
$$
  T_5F^\ell=5^{4\ell-1}\sum_{j=0}^5F^\ell\Big|_{8\ell}\gamma_j.
$$
Then an application of Newton's identity yields
\begin{equation*}
\begin{split}
  \prod_{j=0}^5\left(1-\frac{F\big|_8\gamma_j}{F}y\right)
&=1+\frac{114}{125}y-\frac{6333}{78125}y^2
 +\frac{4}{5^{11}}(8640000t-5177953)y^3 \\
&\qquad+\frac{3}{5^{15}}(8467200000t+1804020097)y^4 \\
&\qquad+\frac{726}{5^{20}}(93744000000t-3501556201)y^5 \\
&\qquad+\frac1{5^{16}}(138240t+14641)^2y^6.
\end{split}
\end{equation*}
Replacing $t$ by $x$ in the last expression, we get the polynomial
$Q(x,y)$. Computing the resultant of the two polynomials
in \eqref{equation: Q} with respect to $y$, we find that the
polynomial $P(x,z)$ is
\begin{equation} \label{equation: P}
\begin{split}
  P(x,z)&=27x^3(x-1)^2(138240x+14641)z^6+7464960x^3(x-1)^2z^5 \\
  &\qquad+1080x^2(x-1)(5760x-5041)z^4+8640x(x-1)(320x-183)z^3 \\
  &\qquad+2160(320x^2-363x+75)z^2+18432(5x-3)z+5120.
\end{split}
\end{equation}
By \cite{Elkies,Errthum}, we know that $t(\tau_0)=-7^4/15^3$ and
consequently, $f(\tau_0)$ is equal to one of the six zeros of the
polynomial $P(-7^4/15^3,z)$ in $z$. To determine which zero is
$f(\tau_0)$, we note that the polynomial $P(x,z)$ can actually be
taken as a defining equation of the Shimura curve $X_6(5)/W_6$ over
$\Q$, where $X_6(5)$ denotes the Shimura curve associated to the
Eichler order $\O\cap(\gamma^{-1}\O\gamma)$ of level $5$ and
$X_6(5)/W_6$ denotes its quotient by the Atkin-Lehner involution $w_2$
and $w_3$. Since $\tau_0$ as a point on $X_6^\ast(5)$ is the unique
CM-point of discriminant $-120$, it is a rational point on
$X_6^\ast(5)$ over $\Q$. In other words, the only rational zero
$2250/6517$ of the polynomial $P(-7^4/15^3,z)$ must be the value of
$f(\tau_0)$. (The other five zeros correspond to the CM-points of
discriminant $-5^2\cdot120$ on $X_6^\ast(5)$.) Thus, the two numbers
$R$ and $S$ in \eqref{equation: RS} in the case of $D=-120$ are
$$
  R=\frac{2^2\cdot7^3\cdot19}{15^{15/4}}, \qquad
  S=\frac{2250}{6517},
$$
and \eqref{equation: theorem 2 3} and \eqref{equation: theorem 2 4}
become
$$
  \frac{2^2\cdot7^3\cdot19}{15^{5/4}}\sum_{n=0}^\infty
  \left(8\cdot\frac{15^3}{7^4}n+\frac{2250}{6517}\right)A_n
  \left(-\frac{7^4}{15^3}\right)^n
 =\frac2{\sqrt5}C_1
$$
and
$$
  \frac{2^2\cdot7^3\cdot19}{15^{5/4}}\sum_{n=0}^\infty
  \left(8\cdot\frac{15^3}{7^4}(n+1/2)+\frac{2250}{6517}\right)A_n'
  \left(-\frac{7^4}{15^3}\right)^n
 =\frac{30\sqrt3}{49}C_1^{-1}.
$$
This proves the identities for the discriminant $-120$.

We next consider the case $D=-19$. Choose the optimal embedding
$\phi:\Q(\sqrt{-19})\hookrightarrow B$ relative to
$(\O,\Z[(1+\sqrt{-19})/2])$ to be the one determined by
$\phi(\sqrt{-19})=5I-J+IJ$. Let $\gamma_0=\iota(\phi(\sqrt{-19}))$ and
write it as $\gamma_0=\SM abcd$. Its fixed point
$\tau_0=(-\sqrt3+\sqrt{-19})/(5+\sqrt3)$ in the upper half-plane lies 
in the fundamental domain. The main task is to find the value of
$$
  S'=s'(\tau)^{-1}\frac d{d\tau}\log\frac{G(\tau)}{\wt G(\tau)}
$$
in \eqref{equation: R'S'} at $\tau_0$, where $G(\tau)$ is any nonzero
automorphic form of weight $12$ on $X_6^\ast$, $\wt
G=G\big|_{12}\gamma_0$, and $s(\tau)=1/t(\tau)$. By Lemma \ref{lemma:
  basis}, we know that we can choose $G=t^{1/2}F^{3/2}$. Thus,
\begin{equation} \label{equation: -19 S}
\begin{split}
  S'&=-t'(\tau_0)^{-1}t(\tau_0)^2\left(
    \frac32\frac d{d\tau}\log\frac{F(\tau)}{\wt F(\tau)}
    \Bigg|_{\tau=\tau_0}+\frac{t'(\tau_0)}{2t(\tau_0)}
   -\frac{19}{(c\tau_0+d)^2}\frac{t'(\gamma_0\tau_0)}{t(\gamma_0\tau_0)}
    \right) \\
 &=-\frac32t'(\tau_0)^{-1}t(\tau_0)^2\frac d{d\tau}
    \log\frac{F(\tau)}{\wt F(\tau)}\Big|_{\tau=\tau_0}-t(\tau_0),
\end{split}
\end{equation}
where $\wt F=F\big|_8\gamma_0$.

Now observe that
$$
  \gamma_1:=\iota\left(\phi\left(\frac{1+\sqrt{-19}}2\right)\right)
  =\frac12(I_2+\gamma_0)=\frac12\M{1+a}bc{1+d}.
$$
has determinant $5$ and can be taken to be a coset representative for
$\Gamma^\ast\backslash\Gamma_5^\ast$ defining the Hecke operator
$T_5$. Therefore, letting $\wt F_1(\tau)=F(\tau)\big|_8\gamma_1$ and
$$
  f_1(\tau)=t'(\tau)\frac d{d\tau}\log\frac{F(\tau)}{\wt F_1(\tau)},
$$
we have $P(t,f_1)=0$, where $P(x,z)$ is the polynomial in
\eqref{equation: P}. By \cite{Elkies,Errthum}, we have
$t(\tau_0)=-3^7/2^{10}$ and $f_1(\tau_0)$ is equal to one of the six
zeros of the polynomial $P(-3^7/2^{10},z)$ in $z$. To see which zero
is equal to $f_1(\tau_0)$, we note that $\tau_0$, as a point on
$X_6^\ast(5)$, is one of two CM-points of discriminant $-19$ and the
these two CM-points of discriminant $-19$ are rational over
$\Q(\sqrt{-19})$. Therefore, $f_1(\tau_0)$ is equal to one of
$512(19\pm\sqrt{-19})/60021$. Now
$$
  \frac d{d\tau}\log\frac{F(\tau)}{\wt F_1(\tau)}
 =\frac{F'(\tau)}{F(\tau)}-\frac{20}{(c\tau+d+1)^2}
  \frac{F'(\gamma_1\tau)}{F(\gamma_1\tau)}+\frac{8c}{c\tau+d+1}.
$$
Evaluating at $\tau_0$, we get
$$
  \frac d{d\tau}\log\frac{F(\tau)}{\wt F_1(\tau)}\Big|_{\tau=\tau_0}
 =\frac{2\sqrt{-19}}{1+\sqrt{-19}}\frac{F'(\tau_0)}{F(\tau_0)}
 +\frac{8c}{1+\sqrt{-19}}.
$$
Comparing with \eqref{equation: main 1}, we find
$$
  \frac d{d\tau}\log\frac{F(\tau)}{\wt F(\tau)}\Big|_{\tau=\tau_0}
 =\frac{1+\sqrt{-19}}{\sqrt{-19}}\frac d{d\tau}\log
  \frac{F(\tau)}{\wt F_1(\tau)}\Big|_{\tau=\tau_0}.
$$
Since $S'$ is necessarily real, we conclude that
$$
  \frac d{d\tau}\log\frac{F(\tau)}{\wt F_1(\tau)}\Big|_{\tau=\tau_0}
 =\frac{512(19+\sqrt{-19})}{60021}, \qquad
  \frac d{d\tau}\log\frac{F(\tau)}{\wt F(\tau)}\Big|_{\tau=\tau_0}
 =\frac{10240}{60021}.
$$
Substituting this into \eqref{equation: -19 S}, we obtain
$S'=15309/15808$ and \eqref{equation: theorem 2 7} and
\eqref{equation: theorem 2 8} become
$$
  \frac{2^{25/3}\cdot13\cdot19^{1/2}}{3^{28/3}}\sum_{n=0}^\infty
  \left(12\cdot\frac{3^7}{2^{10}}n+\frac{15309}{15808}\right)B_n
  \left(-\frac{2^{10}}{3^7}\right)^n=\frac{4\sqrt2}{\sqrt{19}}C_2
$$
and
$$
  \frac{2^{25/3}\cdot13\cdot19^{1/2}}{3^{28/3}}\sum_{n=0}^\infty
  \left(12\cdot\frac{3^7}{2^{10}}(n+1/3)+\frac{15309}{15808}\right)B_n'
  \left(-\frac{2^{10}}{3^7}\right)^n=\frac{3^{7/3}}{2^{5/6}\cdot19^{1/2}}C_2^{-1}.
$$
This proves the identities for the discriminant $-19$.
\end{section}
\bigskip

\centerline{\sc Appendix. List of the auxiliary polynomials in the
  proof of Theorem \ref{theorem: identities}}
\bigskip

Here we list the polynomials used to evaluate the quantities $S$ and
$S'$ in Theorem \ref{theorem: main identities}, corresponding to the
Hecke operators $T_5$, $T_7$, $T_{11}$, and $T_{13}$. We have also
computed the polynomials corresponding to $T_{17}$ and $T_{19}$, but
they are too complicated to be displayed here.

\begin{equation*}
\begin{split}
  P_5(x,z)&=27x^3(x-1)^2(138240x+14641)z^6+7464960x^3(x-1)^2z^5 \\
  &\qquad+1080x^2(x-1)(5760x-5041)z^4+8640x(x-1)(320x-183)z^3 \\
  &\qquad+2160(320x^2-363x+75)z^2+18432(5x-3)z+5120.
\end{split}
\end{equation*}
\begin{equation*}
\begin{split}
  P_7(x,z)&=22235661x^6(x-1)^3(3024000000x-4097152081)z^8\\
  &\qquad+134481277728000000x^6(x-1)^3z^7 \\
  &\qquad+31129925400x^5(x-1)^2(3600000x-1210687)z^6 \\
  &\qquad+1245197016000x^4(x-1)^2(40000x-3471)z^5\\
  &\qquad+648270000x^3(x-1)(19208000x^2-20441799x+253587)z^4\\
  &\qquad+484041600000x^3(x-1)(3430x-4083)z^3\\
  &\qquad+196000000x^2(470596x^2-1267231x+809757)z^2\\
  &\qquad-77760000000x(70x-79)z+97200000000
\end{split}
\end{equation*}
\begin{equation*}
\begin{split}
  P_{11}(x,z)&=
  3^3\cdot11^{11}x^9(x-1)^5(55427328000000000000x^2 \\
  &\qquad  -49446923464224000000x+16546678259573027281)z^{12} \\
  &+146426514344294976000000x^9(x-1)^5(11664000000x-5202748681)z^{11}
  \\
  &+3081366042598800x^8(x-1)^4(1016167680000000000x^2 \\
  &\qquad-1049726627175600000x+115543881567821837)z^{10} \\
  &+616273208519760000x^7(x-1)^4(5645376000000000x^2 \\
  &\qquad-3677542804408000x+127465483643709)z^9 \\
  &+115753795740000x^6(x-1)^3(22541986368000000000x^3 \\
  &\qquad-29303542018872968000x^2+8367736299973548993x \\
  &\qquad-48085474493471259)z^8 \\
  &+67908893500800000x^6(x-1)^3(20492714880000000x^2\\
  &\qquad-20033338665168220x+3718941297294093)z^7 \\
  &+7086244000000x^5(x-1)^2(76372249814784000000x^3 \\
  &\qquad-128575618679598352328x^2+57998197119448630451x \\
  &\qquad-5663679419932760835)z^6 \\
  &+7653143520000000x^4(x-1)^2(20204298892800000x^3 \\
  &\qquad-28634937172769840x^2+9248623134907827x \\
  &\qquad-457073405776467)z^5 \\
  &+4392300000000x^3(x-1)(7334160498086400000x^4 \\
  &\qquad-15976959347457604672x^3+10447013111651394091x^2 \\
  &\qquad-1816427912401084950x+28962670375592235)z^4 \\
  &+556592256000000000x^3(x-1)(8574355240000x^3 \\
  &\qquad-16864516118791x^2+9292154719289x-906852132306)z^3 \\
  &+87120000000000x^2(5477984075731200x^4-15242678431737988x^3 \\
  &\qquad+14384207946457063x^2-4779637835923170x+160347688352127)z^2
  \\
  &+2^{18}\cdot3^2\cdot5^{11}x(251074270137680x^4-659218490328772x^3 \\
  &\qquad+581453582393871x^2-172920969584946x-313911658089)z \\
  &+2^{16}\cdot5^{12}(50214854027536x^4-125349434637768x^3 \\
  &\qquad+105702451736409x^2-30494063413242x+3^6\cdot31^4).
\end{split}
\end{equation*}
\begin{equation*} \small
\begin{split}
  P_{13}(x,z)&=220795952705752437x^9(x-1)^7(4298005400832000000000000x^2\\
  &\qquad+13214513377973804832000000x+86160445679273570730609121)z^{14} \\
  &+4959960281582022744768000000x^9(x-1)^7 \\
  &\qquad(765314352000000x+1176505820688551)z^{13} \\
  &+2296277908139825344800x^8(x-1)^6(3030644833920000000000x^2 \\
  &\qquad-783773604230560200000x-5630164854296621557201)z^{12} \\
  &+229627790813982534480000x^7(x-1)^6(33673831488000000000x^2 \\
  &\qquad-17579212510185032000x-8618813715246332883)z^{11} \\
  &+2009976811158420000x^6(x-1)^5(2885274903386304000000000x^3 \\
  &\qquad-4103184318785684647112000x^2+2044007448506382902634455x\\
  &\qquad+30157407046005622302021)z^{10} \\
  &+1393583922403171200000x^6(x-1)^5(2219442233374080000000x^2\\
  &\qquad-2426222895096718949060x+787891279616791949667)z^9 \\
  &+88098917868000000x^5(x-1)^4(13653120842823990528000000x^3\\
  &\qquad-22388427898915331670097928x^2+9169755921417583303317311x\\
  &\qquad-1335510272388880812298095)z^8 \\
  &+2439662340960000000x^4(x-1)^4(140865532505326886400000x^3 \\
  &\qquad-160884774367727917569920x^2+30658539719491288264803x\\
  &\qquad-2202805319008915077183)z^7 \\
  &+169420995900000000x^3(x-1)^3(422596597515980659200000x^4\\
  &\qquad-702320292170391172660096x^3+322825169771624225610931x^2\\
  &\qquad-17212195753867123465230x+521150155117197321483)z^6 \\
  &+108429437376000000000x^3(x-1)^3(97823286462032560000x^3\\
  &\qquad-122470413622723541518x^2+40960825408914015813x \\
  &\qquad+487853844789766851)z^5\\
  &+30845880000000000x^2(x-1)^2(34386841657133685491200x^4 \\
  &\qquad-67942078041333868279336x^3+39789568931225592998025x^2 \\
  &\qquad-7309981850726385482262x-213332931298137524427)z^4 \\
  &+1645113600000000000x(x-1)^2(39075956428561006240x^4 \\
  &\qquad-77239218331736761300x^3+33148967185029083505x^2\\
  &\qquad-3997021281851275353x-65453273808200820)z^3 \\
  &+676000000000000(x-1)(2641534654570724021824x^5 \\
  &\qquad-10403175406132029675440x^4+9904647406736445306112x^3\\
  &\qquad-2232299734677813245211x^2+121697198509392431190x \\
  &\qquad-992100687686581707)z^2 \\
  &-622080000000000000(x-1)(202316342455567340x^3-249951695360269903x^2\\
  &\qquad+27606248088252470x-201317498982675)z \\
  &+1555200000000000000(1445116731825481x^2-1135066723251890x \\
  &\qquad+7456203666025).
\end{split}
\end{equation*}

\end{document}